\documentclass[12pt, oneside, a4paper]{article}

\usepackage{amsmath}
\usepackage{amsfonts}
\usepackage{amssymb}
\usepackage{amsthm,mathrsfs}
\usepackage{enumerate}
\usepackage{graphicx}
\newtheorem{theorem}{Theorem}[section]

\newtheorem{corollary}[theorem]{Corollary}

\theoremstyle{definition}

\title{\textbf{Cayley graphs on symmetric groups generated by $n$-cycles are hyperenergetic}}
\author{Mahdi Ebrahimi\footnote{ m.ebrahimi.math@ipm.ir}
 \\
 {\small\em  School of Mathematics, Institute for Research in Fundamental Sciences (IPM)},\\{\small\em P.O. Box: 19395--5746, Tehran, Iran},\\{\small\em ORCID ID: 0000-0001-9789-7376}\\% Declarations of interest: none
\\
}
\date{}

\begin{document}

\maketitle

%\fontsize{22}{23}\selectfont
%\baselineskip=12mm

\begin{abstract}
Let $\Gamma$ be a simple graph with $n$ vertices. The energy of $\Gamma$, denoted by $\mathcal{E}(\Gamma)$, is defined as the sum of the absolute values of the eigenvalues of $\Gamma$.
The graph $\Gamma$ is said to be hyperenergetic if $\mathcal{E}(\Gamma)> 2n-2$. For the graph $\Gamma$, the multiplicity of the eigenvalue $0$, denoted by $\eta(\Gamma)$, is called the nullity of $\Gamma$.
  In this paper, we show that for every positive integer $n\geq 4$, the Cayley graph $\Gamma_n$ on the symmetric group $\mathrm{Sym}(n)$ generated by $n$-cycles is an integral hyperenergetic graph with $\mathcal{E}(\Gamma_n)=2^{n-1}(n-1)!$ and $\eta(\Gamma_n)=n!-\binom{2n-2}{n-1}$.

   \end{abstract}

\section{Introduction}
$\noindent$In this paper, all groups and graphs are assumed to be finite.
 Let $\Gamma$ be a simple graph with vertex set $\{\nu_1,\nu_2,\dots, \nu_n\}$.
 The \textit{adjacency matrix} of $\Gamma$, denoted by $A(\Gamma)$, is the $n\times n$ matrix such that the $(i,j)$-entry is $1$ if $\nu_i$ and $\nu_j$ are adjacent, and is $0$ otherwise.
  The \textit{eigenvalues} of $\Gamma$ are the eigenvalues of its adjacency matrix $A(\Gamma)$.

Suppose $\{\lambda_1, \lambda_2, \dots, \lambda_n\}$ is the set of all eigenvalues of $\Gamma$.
 The \textit{energy} of the graph $\Gamma$ is defined as $\mathcal{E}(\Gamma):=\sum_{i=1}^n|\lambda_i|$.
  This concept was introduced by Gutman \cite{182}.
  The total $\pi$-electron energy and various "resonance energies" derived from it, play an essential role in the theory of conjugated molecules, for more details see \cite{8}.
   In recent decades, the energy of a graph has been much studied by researchers, for instance see, \cite{Ak,7,GZ,oboudi,zh}.

The simple graph $\Gamma$ with $n$ vertices is said to be \textit{hyperenergetic} if its energy exceeds the energy of the complete graph $K_n$; that is, if $\mathcal{E}(\Gamma)>2n-2$. Otherwise, $\Gamma$ is called \textit{non-hyperenergetic}. Hyperenergetic graphs were first introduced by Gutman \cite{Hyper}. It has been shown \cite[Corollary 7.8]{9} that for every positive integer $n\geq 8$, there exists a hyperenergetic graph of order $n$. In recent years, a very large number of papers on hyperenergetic graphs has been published (see \cite{kn,9,13,Hyper,10,15}). In this paper, we wish to prove that Cayley graphs on symmetric groups generated by $n$-cycles are hyperenergetic.

For the graph $\Gamma$, the multiplicity of the eigenvalue $0$, denoted by $\eta(\Gamma)$, is called the \textit{nullity} of $\Gamma$. The study of graph nullity was motivated by its applications to chemistry. In 1950, Longuet-Higgings \cite{LH} showed that if the nullity of the molecular graph of an alternant hydrocarbon is non-zero, then the molecule is unstable. In 1957, Collatz and sinogowit \cite{CS} proposed the problem of classifying all graphs with nullity zero. Afterwards, this problem and relevant questions have been extensively studied, for instance, see \cite{cheng,Go,Ru,Ta,Wa}.

 Let $G$ be a finite group and $S$ be
 an inverse closed subset of $G$ with $1 \notin S$.
  The \textit{Cayley graph} $\mathrm{Cay}(G,S)$ is the graph which has the elements of
   $G$ as its vertices and two vertices $u,\nu \in G$
    are joined by an edge if and only if $\nu=au$, for some $a\in S$. A Cayley graph $\mathrm{Cay}(G,S)$ is called \textit{normal} if $S$ is closed under conjugation with elements of $G$. Also $\mathrm{Cay}(G,S)$ is called \textit{integral} if its eigenvalues are all integers.

 \begin{theorem}\label{main}
 For a positive integer $n\geq 4$, let $G$ be the symmetric group $\mathrm{Sym}(n)$ on $n$ letters, and $Z_n$ be the set of all $n$-cycles of $G$. Then the Cayley graph $\Gamma_n:=\mathrm{Cay}(G,Z_n)$ is an integral hyperenergetic graph with  $\mathcal{E}(\Gamma_n)=2^{n-1}(n-1)!$ and $\eta(\Gamma_n)=n!-\binom{2n-2}{n-1}$.
  \end{theorem}

To determine the significance of Theorem \ref{main}, we use it to prove the following interesting result.
\begin{corollary}\label{binom}
For every positive integer $n$,
$$a_n:=\sum_{k=0}^{n-1}(\binom{2n}{n}-2\binom{2n}{k})$$
is a non-negative integer.
\end{corollary}

 %%%%%%%%%%%%%%%%%%%%%%%%%%%%%%%%%%%%%%
\section{Cayley graphs on symmetric groups generated by $n$-cycles are hyperenergetic}
$\noindent$We first state well-known results on the character theory of the symmetric groups; for a complete account, see \cite{GA}.
A partition of a positive integer $n$ is a tuple $\alpha =(\alpha_{1},\alpha_{2},\dots,\alpha_{r})$ of positive integers $\alpha_{1}\geq \alpha_{2}\geq \dots\geq \alpha_{r}$ such that $\alpha_{1}+\alpha_{2}+ \dots+\alpha_{r}=n$.
  The integers $\alpha_i$'s are called the parts of $\alpha$.
    To indicate that $\alpha$ is a partition of $n$, we write $\alpha \vdash n$.  For $i = 1,\dots , n$, if $t_i$ is the number of parts of $\alpha$ equal to $i$, then we can also write
$\alpha=(r^{t_r},\dots,1^{t_1})$.  Usually $i^{t_i}$ is left out if $t_i = 0$.
If $\alpha =(\alpha_{1},\alpha_{2},\dots,\alpha_{r})$ is a partition of $n$ then the Young diagram  $[\alpha]$ of $\alpha$ consists of $n$ boxes placed into $r$ rows, where the $i$-th row has $\alpha_{i}$ boxes.  The box in the $i$-th row and $j$-th column of $[\alpha]$ is called the $(i,j)$ node of $[\alpha]$.
For each $i$, denote $\alpha^{\top}_i$
the number of parts of $\alpha$ which are bigger than or equal to $i$.  The partition $\alpha^{\top} =(\alpha^{\top}_1,\alpha^{\top}_2,\dots,\alpha^{\top}_s)$ is called
the conjugate partition associated with $\alpha$. If $\alpha=\alpha^\top$, then $\alpha$ is called a self-conjugate partition.

 If $(i, j)$ is a node of $[\alpha]$ we denote by $H^{\alpha}_{i,j}$  the $(i, j)$-hook
of $\alpha$ which is the set of nodes of $[\alpha]$ of the form $(i ,j')$ for some $ j'\geq j$  or $(i', j)$
for some $i'\geq i$.  The hook-length $h^{\alpha}_{i,j}$
 of the $(i,j)$-node is equal to the number of nodes in $H^{\alpha}_{i,j}$. The set of nodes $(l, k)$'s of $\alpha$ with $l\geq i$ and $k\geq j$, satisfying $(l + 1, k + 1)\notin [\alpha]$ is called the (i,j)-rim of $[\alpha]$ and is denoted by $R^{\alpha}_{i,j}$.
   For a partition $\alpha$ of $n$ and $k\in \mathbb{N}$,  we also set $ I_k^{\alpha}:=\{{(i,j)}\,|\,h^\alpha_{i,j}=k\}$.

In the symmetric group $\mathrm{Sym}(n)$, each conjugacy class of $\mathrm{Sym}(n)$ corresponds
naturally to the partitions of $n$ associated to the cycle structure of that class.   The value of the irreducible character $\chi^\alpha$, labeled by the partition $\alpha$, evaluated at the conjugacy class corresponding to a partition $\beta$ can be calculated recursively by the well known Murnaghan-Nakayama formula \cite{GA}.  Precisely,
if $\alpha$ is a partition of $n=k+m$  and $\beta\in \mathrm{Sym}(n)$ contains a $k$-cycle in its cycle structure and $\rho\in \mathrm{Sym}(m)$ is of cycle type deleting $k$-cycle out of $\beta $, then
 $$\chi^{\alpha}(\beta)=\sum_{(i,j )\in I^\alpha_k}(-1)^{l^{\alpha}_{i,j}}\chi^{\alpha\backslash R^{\alpha}_{i,j} }(\rho),$$
where $l^{\alpha}_{i,j}:=\alpha_j^\top-i$ is the leg length of the hook $H^{\alpha}_{i,j}$ and $\alpha\backslash R^{\alpha}_{i,j}$ simply denotes the partition associated to the Young diagram $[\alpha]\backslash R^{\alpha}_{i,j}$. Note that when $I^\alpha_k=\emptyset$, then $\chi^{\alpha}(\beta)=0$.

It is well known that the eigenvalues of a normal Cayley graph  $\mathrm{Cay}(G,S)$ can be expressed in terms of the irreducible characters of the group $G$ \cite[p.235]{eigen}.

\begin{theorem}\label{eigen}(\cite{2,6})
The eigenvalues of a normal Cayley graph $\mathrm{Cay}(G,S)$
are given by $\eta_\chi=\sum_{a\in S}\chi(a)/\chi(1)$ where
 $\chi$ ranges over the set of all complex irreducible characters $\mathrm{Irr}(G)$ of $G$. Moreover,
  the multiplicity of $\eta_{\chi}$ is $\chi(1)^2$.
\end{theorem}
Now we are ready to prove our main result.\\
\textit{Proof of Theorem \ref{main}:}  Applying \cite[Theorem 1]{ko}, the graph $\Gamma_n$ is integral. We now consider some particular partitions of $n$. Assume that $\alpha_0:=(1^n)$, for every integer $1\leq m\leq n-2$, $\alpha_m:=(m+1,1^{n-m-1})$ and $\alpha_{n-1}:=(n)$. It is clear that the conjugacy class $Z_n$ corresponds to the partition $\beta:=(n)$. Let $\alpha\vdash n$. If $\alpha \neq \alpha_m$, for every $0\leq m\leq n-1$, then using  Murnaghan-Nakayama formula, we deduce that $\chi^{\alpha}(\beta)=0$. Now suppose $\alpha=\alpha_m$, for some $0\leq m\leq n-1$. Then again applying Murnaghan-Nakayama formula, we have $\chi^{\alpha}(\beta)=(-1)^{n-m-1}$. Also it is easy to see that $|Z_n|=(n-1)!$. Let $A:=\{\chi^{\alpha_m}|\,0\leq m\leq n-1\}$. Then by Theorem \ref{eigen},
\begin{align}
\mathcal{E}(\Gamma_n)&=\sum_{\chi \in \mathrm{Irr}(G)}\chi(1)^2|\eta_\chi|\nonumber\\
&=\sum_{\chi \in \mathrm{Irr}(G)-A}\chi(1)^2|\eta_\chi|+\sum_{m=0}^{n-1}\chi^{\alpha_m}(1)^2|\eta_{\chi^{\alpha_m}}|\nonumber\\
&=0+\sum_{m=0}^{n-1}\chi^{\alpha_m}(1)^2|\eta_{\chi^{\alpha_m}}|\nonumber\\
&=\sum_{m=0}^{n-1}\chi^{\alpha_m}(1)^2|\eta_{\chi^{\alpha_m}}|.\nonumber
\end{align}
Now suppose $m\leq n-1$ is a non-negative integer.
Then using Murnaghan-Nakayama formula, $\chi^{\alpha_m}(1)=\binom{n-1}{m}$.
 Hence $\eta_{\chi^{\alpha_m}}=(-1)^{n-m-1}(n-1)!/\binom{n-1}{m}$. Thus
 \begin{align}
\mathcal{E}(\Gamma_n)&=(n-1)!\sum_{m=0}^{n-1}\binom{n-1}{m}\nonumber\\
&=2^{n-1}(n-1)!.\nonumber
\end{align}
Hence $\Gamma_n$ is hyperenergetic. Also using Theorem \ref{eigen}, it is clear that
\begin{align}
\eta(\Gamma_n)&=|G|-\sum_{m=0}^{n-1}\chi^{\alpha_m}(1)^2\nonumber\\
&=n!-\sum_{m=0}^{n-1}\binom{n-1}{m}^2.\nonumber
\end{align}
Therefore by Vandermonde's identity, we deduce that
$$\eta(\Gamma_n)=n!-\binom{2n-2}{n-1}.\qed$$

We end this section by the proof of Corollary \ref{binom}.\\
\textit{Proof of Corollary \ref{binom}:}
If $n=1$ or $2$, then it is easy to see that $a_n$ is non-negative. Thus we can assume that $n\geq 3$. Applying Theorem \ref{main} for the symmetric group $\mathrm{Sym}(n+1)$ and \cite[Theorem 4.5]{Li},
$$2^n n!\leq \sqrt{(n+1)!n!\binom{2n}{n}}.$$
Thus
$$2^{2n}\leq (n+1)\binom{2n}{n}.$$
Hence as
\begin{align}
(n+1)\binom{2n}{n}-2^{2n}&=(n+1)\binom{2n}{n}-\sum_{k=0}^{2n}\binom{2n}{k}\nonumber\\
&=(n+1)\binom{2n}{n}-\sum_{k=0}^{n-1}2\binom{2n}{k}-\binom{2n}{n}\nonumber\\
&=\sum_{k=0}^{n-1}(\binom{2n}{n}-2\binom{2n}{k})\nonumber\\
&=a_n,\nonumber
\end{align}
We deduce that $a_n$ is a non-negative integer.\qed
%%%%%%%%%%%%%%%%%%%%%%%%%%%%%%%%%%%%%%%%%%%%%%

\section*{Acknowledgements}
Funding: This research was supported in part
by a grant  from School of Mathematics, Institute for Research in Fundamental Sciences (IPM).
%%%%%%%%%%%%%%%%%%%%%%%%%%%%

%%%%%%%%%%%%%%%%%%%%%%%%%%%%%

\end{document}